\documentclass[a4paper,12pt]{article}
\usepackage[dvips]{graphicx}
\DeclareGraphicsExtensions{.eps,.ps,.eps.gz,.ps.gz}
\newcount\vgl
\newcount\vgla
\newcount\vglb
\catcode`@=11
\def\eqalijn#1{\null\,\vcenter{\openup1\jot \m@th
       \ialign{\strut\hfil$\displaystyle{##}$&$\displaystyle{{}##}$\hfil
               &$\displaystyle{{}##}$\hfil&$\displaystyle{{}##}$\hfil
               &$\displaystyle{{}##}$\hfil&$\displaystyle{{}##}$\hfil
               \crcr#1\crcr}}\,}
\def\dod#1{\global\advance\vgl by 1

      \expandafter\expandafter\global\setbox\vgl=\hbox{$\displaystyle{#1}$}}
\def\doa#1{\setbox0=\hbox{$\displaystyle{{}#1}$}}
\def\dob{\ifdim \ht0>\expandafter\ht\vgl \expandafter\ht\vgl=\ht0\fi
                \ifdim \dp0>\expandafter\dp\vgl \expandafter\dp\vgl=\dp0\fi}
\def\initia{\vgl=0\vglb=0}
\def\eqalijna#1{\null\,\vcenter{\openup1\jot\m@th\ialign{
  \strut{}\hfil\dod{##}\copy\vgl%
       &\doa{##}\copy0\hfil\dob
 &\doa{##}\copy0\hfil\dob
       &\doa{##}\copy0\hfil\dob
       &\doa{##}\copy0\hfil\dob
       &\doa{##}\copy0\hfil\dob\crcr#1\crcr}}\,}
\def\eqalijnb#1{\null\vcenter{\openup1\jot\m@th
        \ialign{\strut\global\advance\vglb by 1$\vphantom{\copy\vglb}${##}&
                \hfil$\smash{##}$\crcr#1\crcr}}}
\catcode`@=12 \initia
\def\latexeqno#1{\refstepcounter{equation}\label{#1}(\theequation)}
\def\latexeqno#1{\refstepcounter{equation}\label{#1}(\theequation)}
\def \a {\alpha}
\def \i {{\cal{J}}}
\def \k {{\cal{K}}}
\def \dn {{\delta_n}}
\def \d {\delta}
\def \cd {c_\d}

\def \calR  {{\cal{R}}^\mu}

\def \P  {{\cal{P}}}
\def \Pmu  {{\cal{P}}^\mu}
\def \psimu  {\psi^\mu}
\def \calQ  {{\cal{Q}}^\mu}
\def \QT {{Q_T}}
\def \Mmu {M^\mu}

\def \e {\varepsilon}

\def \al {\alpha}

\def \fmu {f^{\mu}}
\def \kw {k_w}
\def \ka {k_a}

\def \su {\underline s}
\def \so {\overline s}
\def \umu {u^{\mu}}
\def \pmu {p^{\mu}}

\def \O {{\Omega}}
\def \pc {p_c}
\def \N {N}
\def \R {{\mathrm {I\mkern-5.5mu R\mkern1mu}}}

\newcommand\bes{\begin{eqnarray}}
\newcommand\ees{\end{eqnarray}}
\newcommand\bess{\begin{eqnarray*}}
\newcommand\eess{\end{eqnarray*}}
\newcommand{\Int} {\displaystyle \int}
\newcommand{\IntO}{\displaystyle \Int _{\Omega}}

\newcommand{\mint}{\displaystyle {\Int \kern -0.961em -}}
\newcommand{\mintO}{\displaystyle {\Int \kern -0.961em -}_{\Omega}}
\newcommand{\mintG}{\displaystyle {\Int \kern -0.961em -}_{\Gamma}}

\newcommand{\Frac}[2] {\frac{\textstyle #1} {\textstyle #2}}
\renewcommand{\theequation}{\thesection.\arabic{equation}}
\newcommand{\proof}{{\underline{Proof}}: }
\newtheorem{definition}{Def\hskip 1pt inition }[section]
\newtheorem{lemma}[definition]{Lemma}

\newtheorem{theorem}[definition]{Theorem}
\newtheorem{corollary}[definition]{Corollary}
\begin{document}
\title{Singular limit of a two-phase flow problem in porous medium as the air viscosity tends to zero\thanks{This work was supported by the
GNR MoMaS (PACEN/CNRS, ANDRA, BRGM, CEA, EdF, IRSN), France.}}
\date{}
\maketitle{}
\begin{center}\author{R. Eymard \footnote{Universit\'e Paris-Est Marne-La-Vall\'ee, 5 bd Descartes, Champs-sur-Marne,
77454 Marne-la-Vall\'ee Cedex 2,  France}, \and M. Henry \footnote{CMI Universit\'e de Provence, 39 rue Fr\'ed\'eric
Joliot-Curie 13453 Marseille cedex 13, France} \and and D. Hilhorst \footnote{CNRS and Laboratoire de
Math\'ematiques, Universit\'e de Paris-Sud 11, F-91405 Orsay Cedex, France}
 }
\end{center}

\begin{abstract}
In this paper we consider a two-phase flow problem in porous media and study its singular limit as the viscosity of
the air tends to zero; more precisely, we prove the convergence of subsequences to solutions of a generalized
Richards model.
\end{abstract}

\section{Introduction}\label{I}
Hydrologists have studied air-water flow in soils, mainly using the so-called Richards approximation. At least two
hypotheses are physically required for this model to be applicable: the water pressure in the saturated region must
be larger than the atmospheric pressure and all the unsaturated regions must have a boundary connected to the
surface. However, in many situations, these hypotheses are not satisfied and a more general two-phase flow model
must be considered. This work explores the limit of this general model as the viscosity of the air tends to zero,
which is one of the hypotheses required in the Richards model. To that purpose we prove the existence of a weak
solution of the two-phase flow problem and prove estimates which are uniform in the air viscosity. In this paper, we
assume that the air and water phases are incompressible and immiscible. The geometric domain is supposed to be
horizontal, homogeneous and isotropic. Our starting point is the following two-phase flow model, which one can
deduce from Darcy's law
$$(\cal{TP})\left\{\begin{array}{ll}
&u_t -div(k_w(u)\nabla(p)) = s_w\nonumber\\
&(1-u)_t - div(\Frac{1}{\mu} k_a(u)\nabla(p + p_c(u))) = s_a, \nonumber
\end{array}\right.$$
where $u$ and $p$ are respectively the saturation and the pressure of the water phase, $k_w$ and $k_a$ are
respectively the relative permeabilities of the water and the air phase, $\mu$ is the ratio between the viscosity of
the air phase and that of the water phase, $p_c$ is the capillary pressure, $s_w$ is an internal source term for the
water phase and $s_a$ is an internal source term for the air phase; these source terms are used to represent
exchanges with the outside. We suppose in particular that the physical functions $k_w$, $k_a$ and $p_c$ only depend
on the saturation $u$ of the water phase, and that $ k_w(1) = k_a(0) = 1$. The aim of this paper is the study of the
limit of the two-phase flow problem as $\mu \downarrow 0$.\\
$\\ $
The classical Richards model as formulated by the
engineers is given by
$$({\cal{R}})\left\{\begin{array}{ll}
&u_t -div(k_w(u)\nabla p) = s_w\nonumber\\
&u=p_c^{-1}(p_{atm}-p). \nonumber
\end{array}\right.$$
where the properties of capillary pressure $p_c=p_c(u)$ are describes in hypothesis $(H_8)$ below. For the existence
and uniqueness of the solution of Richards model together with suitable initial and boundary conditions as well as
qualitative properties of the solution and methods for numerical approximations we refer to \cite{AL}, \cite{HW},
\cite{P}, \cite{RPK}. In this article, we will show that the singular limit as $\mu\downarrow 0$ of the two phase
flow problem $(\cal{TP})$ has the form
$$({\cal{FBP}})\left\{\begin{array}{ll}
&u_t -div(k_w(u)\nabla p) = s_w\nonumber\\
&u=1\mbox{ or }\nabla(p+p_c(u))=0\mbox{ a.e. in }
\Omega\times(0,T). \nonumber
\end{array}\right.$$
We remark that a solution of $({\cal{R}})$ with $u > 0$ satisfies $({\cal{FBP}})$. \\
$\\ $ This paper is organized as follows. In Section 2 we present a complete mathematical formulation of the
problem, and state the main mathematical results, which include a precise formulation of the singular limit problem.
We give a sequence of regularized problems in Section 3, and prove the existence of a classical solution. In Section
4 we present a priori estimates, which are uniform in an extra regularization parameter $\delta$ and in the air
viscosity $\mu$. In Section 5, we let $\delta\downarrow 0$ and prove that the solution converges to a solution of
the two phase flow problem. We study its limiting behavior as the air viscosity $\mu$ tends to zero in Section 6.
Finally in Section 7 we propose a finite volume algorithm in a one
dimensional context and present a variety of numerical solutions.\\

\section{Mathematical formulation and main results}\label{I}
We consider the two-phase flow problem
$$
\initia \mbox{$(S^{\mu})$} \left\{ \eqalijna{ &u_t =div\bigg(\kw(u)\nabla p\bigg)+\fmu(c)\so-\fmu(u)\su,
~~~~~&\mbox{ in }Q_T,\cr &(1-u)_t = div\bigg(\Frac{1}{\mu}\ka(u)\nabla(p+\pc(u))\bigg)\cr &~~~~~~~~~~~~~~+
(1-\fmu(c))\so-(1-\fmu(u))\su,~~~~&\mbox{ in }Q_T,\cr  &\Int_\O p(x,t)dx=0,&\mbox{ for }t \in (0,T),\cr &\nabla
p.n=0, &\mbox{ on }\partial \Omega \times (0,T),\cr &\nabla (p+\pc(u)).n=0, &\mbox{ on }\partial \Omega \times
(0,T),\cr &u(x,0)=u_0(x),&\mbox{ for }x \in \Omega, } \right. \eqno\eqalijnb{\latexeqno{n1}\cr\nonumber\cr
              \latexeqno{n2} \cr
                \latexeqno{n4} \cr\latexeqno{n3} \cr \latexeqno{n5} \cr \latexeqno{CImu} \cr
                 }
$$
where $T$ is a positive constant, $Q_T:=\Omega \times (0,T)$ and where we suppose that
$$\begin{array}{lll}
&(H_1) &\Omega \mbox{ is a smooth bounded domain of $\R^N$ where the space dimension $N$ is arbitrary},\\
&(H_2) &u_m\in(0,1),\\
&(H_3) &c\in L^\infty(\O\times(0,T)) \mbox{ and }u_m\leq c \leq 1,\\
&(H_4) &u_0\in L^\infty(\O) \mbox{ and }u_m\leq u_0 \leq 1,\\
&(H_5) &\so\in L^2(\Omega),~~\so\geq 0,~~\su\in L^2(\Omega),~~\su\geq 0\mbox{ and }\Int_\O(\so(x)-\su(x))dx=0,\\
&(H_6) &\kw\in C^2([0,1]),~~\kw'\geq 0,~~\kw(0)=0,~~\kw(1)=1\mbox{ and }\kw(u_m)>0,\\
&(H_7) &\ka\in C^2([0,1]),~~\ka'\leq 0,~~\ka(1)=0,~~\ka(0)=1\mbox{ and }\ka(s)>0\mbox{ for all }s\in[0,1),\\
&(H_8) &\pc\in C^0([0,1])\cup C^3([0,1)),~~\pc'<0\mbox{ and }\sup_{s\in[0,1)}(-\ka(s)\pc'(s))<+\infty,\\
&(H_9) &\mu\in(0,1].
\end{array}$$
In this model, $u$ and $p$ are respectively the saturation and the pressure of the water phase, $k_w$ and $k_a$ are
respectively the mobilities of the water phase and the mobility of the non-water phase and $p_c$ is the capillary
pressure. We assume in particular that the permeability functions $k_w$, $k_a$ and the capillary pressure $p_c$ only
depend on the saturation $u$ of the water phase. Here, we suppose that the flow of the water phase in the reservoir
is driven by an injection term $f^\mu(c)\so$ and an extraction term $f^\mu(u) \su$ where $\so$ and $\su$ are given
space dependent functions, $c$ is the saturation of the injected fluid; if $c=1$, only water will be injected, if
$c=0$, only air will be injected, whereas a mixture of water and air will be injected if $0 < c < 1$. The function
$f^\mu$ is the fractional flow of the water phase, namely
\begin{equation}\label{deffmu}
\fmu(s)=\Frac{\kw(s)}{\Mmu(s)}, \mbox{ with
}\Mmu(s)=\kw(s)+\Frac{1}{\mu}\ka(s).
\end{equation}
In particular, we remark that
\begin{equation}\label{fmucroissante}
\fmu(s) \mbox{ is non decreasing. }
\end{equation}
Next we introduce a set of notations, which will be useful in the sequel.
\begin{equation}\label{defg}
g(s)=-\Int_0^s\ka(\tau)\pc'(\tau)d\tau,
\end{equation}
\begin{equation}\label{defzeta}
\zeta(s)=\Int_0^s\sqrt{\ka(\tau)}\pc'(\tau)d\tau,
\end{equation}
\begin{equation}\label{Q}
\calQ(s)=\Int_0^s \fmu(\tau)\pc'(\tau)d\tau,
\end{equation}
and
\begin{equation}\label{R}
\calR(s)=\Int_0^s\Frac{\ka(\tau)}{\ka(\tau)+\mu\kw(\tau)}\pc'(\tau)d\tau,
\end{equation}
for all $s\in[0,1]$. This implies in particular that
\begin{equation}\label{P+Q} \calR(s)+\calQ(s)=\pc(s)-\pc(0), \mbox{ for all }s\in[0,1].
\end{equation}
\begin{definition} The pair $(\umu,\pmu)$ is a weak solution of Problem $(S^\mu)$ if
$$\begin{array}{lll}
&\umu\in L^\infty(\O\times(0,T)),~~~~&\mbox{with~}0\leq \umu\leq 1 \mbox{ in }Q_T,\\
&&\\
&\pmu\in L^2(0,T;H^1(\O)),~~~~&\Int_\O\pmu(x,t)dx=0 \mbox{~for~almost~every~} t \in (0,T),\\
&&\\
&g(\umu)\in L^2(0,T;H^1(\O)),&\\
\end{array}$$
with
\begin{eqnarray}
&\Int_0^T\IntO\umu \varphi_t dxdt=&
\Int_0^T\IntO\kw(\umu)\nabla\pmu.\nabla \varphi dx dt
-\Int_0^T\IntO\bigg(\fmu(c)\so-\fmu(\umu)\su \bigg) \varphi dxdt\nonumber\\
&&-\IntO u_0(x)\varphi(x,0)dx,\label{defsol1}
\end{eqnarray}and
\begin{eqnarray}
&&\Int_0^T\IntO\bigg(1-\umu\bigg)\varphi_tdxdt=\Int_0^T\IntO\Frac{1}{\mu}\ka(\umu)\bigg(\nabla\pmu+\nabla \pc(\umu)\bigg).\nabla \varphi dx dt\nonumber\\
&&-\Int_0^T\IntO\bigg((1-\fmu(c))\so-(1-\fmu(\umu))\su\bigg)
\varphi dx dt-\IntO
\bigg(1-u_0(x)\bigg)\varphi(x,0)dx,~~~~\label{defsol2}
\end{eqnarray}
for all $\varphi$ in ${\cal{C}}:=\{w\in W_2^{2,1}(Q_T), w(.,T)=0 \mbox{ in }\Omega \}.$
\end{definition}
$\\ $ Our first result, which we prove in Section \ref{Existencesolution}, is the following
\begin{theorem}\label{thexistence}
Suppose that the hypotheses $(H_1)-(H_9)$ are satisfied, then there exists a weak solution $(\umu,\pmu)$ of Problem
$(S^\mu)$.
\end{theorem}
Next we define the discontinuous function $\chi$ by $$\chi(s):=\left\{\begin{array}{ll}0&\mbox{ if }s\in[0,1)\\
1&\mbox{ if }s=1,\end{array}\right.$$ as well as the graph
$$H(s):=\left\{\begin{array}{lll}&0&\mbox{ if }s\in[0,1)\\ &[0,1]&\mbox{ if }s=1.\end{array}\right.$$
The main goal of this paper is to prove the following convergence result,
\begin{theorem}\label{thlim}
Suppose that the hypotheses $(H_1)-(H_9)$ are satisfied, then there
exists a subsequence $((u^{\mu_n},p^{\mu_n}))_{n\in N}$ of
weak solutions of Problem $(S^{\mu_n})$ and functions $u$, $p$, $\hat f$ such that
$$\begin{array}{ll}
&u\in L^\infty(Q_T),~~~0\leq u\leq 1\mbox { in }Q_T,\\
&\\
&\hat f\in L^\infty(Q_T),~~~0\leq \hat f\leq 1\mbox { in }Q_T,\\
&\\
&p\in L^2(0,T;H^1(\O)),\\
&\\
&\ka(u)\nabla\pc(u)\in L^2(\O\times(0,T)),\\
\end{array}$$
and
$$\begin{array}{ll}
&(u^{\mu_n})_{n\in N}\mbox{ tends to }u\mbox{ strongly in }L^2(Q_T),\\
&(p^{\mu_n})_{n\in N}\mbox{ tends to }p\mbox{ weakly in }L^2(0,T;H^1(\O)),\\
\end{array}$$
as $\mu_n$ tends to zero and
\begin{eqnarray}
&\Int_0^T\IntO u \varphi_tdxdt =&\Int_0^T\IntO\kw(u)\nabla p.\nabla \varphi dxdt-\Int_0^T\IntO\bigg(\chi(c)\so-\hat f\su \bigg) \varphi dxdt\nonumber\\
&&-\IntO u_0(x)\varphi(x,0)dx,\label{sollim1}
\end{eqnarray}
for all $\varphi \in{\cal{C}}$, where $\hat f(x,t)\in H(u(x,t))$
for $(x,t)\in Q_T$. Moreover we also have that
\begin{equation}\label{lim1}
\Int_0^T\IntO \bigg[\ka(u)\bigg]^2\bigg[\nabla p+\nabla \pc(u)\bigg]^2dxdt=0
\end{equation}
and
\begin{equation}\label{lim2}
\Int_\O p(x,t) dx=0, \mbox{ for almost every }t\in(0,T).
\end{equation}
\end{theorem}
Formally, $u$ satisfies the following limit problem
$$
\left\{
\begin{array}{ll}
u_t =div\bigg(k_w(u)\nabla p \bigg)+ \chi(c)\bar s - \hat
f\underline{s},
~~&\mbox{ in }Q_T, \\
\nabla u.n=0, &\mbox{ on }\partial \Omega \times (0,T),\\
u(x,0)=u_0(x),~~&\mbox{ for }x \in \Omega.
\end {array}
\right.
$$
More precisely the following corollary holds
\begin{corollary}\label{thcoro} Suppose that $u<1$ in ${\cal O}=\cup_{t\in[\tau,T]}\Omega_t$, where $\tau >0$
and $\Omega_t$, for $t\in[\tau,T]$, are smooth subdomains of $\O$
and ${\cal O}$ is a smooth domain of $\Omega\times[\tau,T]$ and that $u=1$ in
$Q_T\setminus\overline{\cal{O}}$ then
$$p(x,t)=-\pc(u(x,t))+constant(t), \mbox{ for all }(x,t)\in {\cal O}$$ and $u$ satisfies
$$
\left\{
\begin{array}{ll}
u_t =-div\bigg(k_w(u)\nabla \pc(u)\bigg)+ \chi(c)\bar s, ~~&\mbox{ in } {\cal O}, \\
\Frac{\partial u}{\partial n}=0,&\mbox{ on }\partial {\cal O}\cap
\bigg(\partial\Omega\times(0,T)\bigg),\\
u=1, &\mbox{ elsewhere on }\partial {\cal O},\\
u(x,0)=u_0(x),&\mbox{ for }x \in \Omega.
\end {array}
\right.
$$
\end{corollary}
$\\ $ Finally we remark that another form of the limit problem involves a parabolic equation, which is close to the
standard Richards equation. Indeed if we set $\phi(s):=\pc(0)-\pc(s)$ and denote by $\beta$ the inverse function of
$\phi$, the function $v:=\phi(u)$ is a weak solution of the problem
$$
\left\{
\begin{array}{ll}
\beta(v)_t =div\bigg(k_w(\beta(v))\nabla v \bigg)+ \chi(c)\bar s -
\hat f\underline{s},
~~&\mbox{ in }Q_T, \\
\nabla \beta(v).n=0, &\mbox{ on }\partial \Omega \times (0,T),\\
\beta(v)(x,0)=u_0(x),~~&\mbox{ for }x \in \Omega,
\end {array}
\right.
$$
with $\hat f\in H(\beta(v))$.
\setcounter{equation}{0}
\section{Existence of a solution of an approximate problem $(S^{\mu}_{\d})$ of Problem
$(S^{\mu})$}\label{Existencesolution}
Let $\d$ be an arbitrary positive constant. In order to prove the
existence of a solution of Problem $(S^{\mu})$ we introduce a
sequence of regularized problems $(S^{\mu}_{\d})$, namely
$$
\initia \mbox{$(S^{\mu}_{\d})$} \left\{ \eqalijna{ &u_t~~ =
div\bigg(\kw(u)\nabla p\bigg)+\fmu(\cd)\so_\d\cr
&~~~~~~~~-\fmu(u)\bigg(\su_\d+\mint_\O(\so_\d-\su_\d )dx\bigg),
~~&\mbox{ in }\Omega \times (0,T),\cr &(1-u)_t
=div\bigg(\Frac{1}{\mu}\ka(u)\nabla(p+\pc(u))\bigg)+\bigg(1-\fmu(\cd)\bigg)\so_\d~~\cr
&~~~~~~~~~~~~~~-\bigg(1-\fmu(u)\bigg)\bigg(\su_\d+\mint_\O(\so_\d-\su_\d
)dx\bigg), &\mbox{ in }\Omega \times (0,T),\cr &\Int_\O
p(x,t)dx=0,&\mbox{ for }t\in(0,T),\cr
 &\nabla p.n=0, &\mbox{ on }\partial \Omega
\times (0,T),\cr &\nabla (p +\pc(u)).n=0, &\mbox{ on }\partial
\Omega \times (0,T),\cr &u(x,0)=u_{0}^\d(x),& \mbox{ for }x \in
\Omega,\cr } \right. \eqno\eqalijnb{\nonumber\cr\latexeqno{1d}\cr
\nonumber\cr
              \latexeqno{2d} \cr
               \latexeqno{3d} \cr
                \latexeqno{4d} \cr   \latexeqno{5d} \cr\latexeqno{CId} \cr
                 }
$$
where  $u_{0}^\d$, $c_\d$, $\so_\d$ and $\su_\d$ are smooth functions such that $u_{0}^\d$ tends to $u_0$ in
$L^2(\Omega)$ and $c_\d$, $\so_\d$ and $\su_\d$ tend respectively to $c$, $\so$ and $\su$ in $L^2(Q_T)$, as
$\d\downarrow 0$. In particular we suppose that there exists a positive constant $C$ such that
\begin{equation}\label{sdeltaborne}
\su_\d\geq 0,~~\so_\d\geq 0 \mbox{ and }\Int_\O \su_\d^2+\Int_\O \so_\d^2 \leq C.
\end{equation}
Moreover we suppose that $u_{0}^\d$, $c_\d$ satisfy
\begin{eqnarray}
&0<u_m\leq u_0^\d \leq 1-\d<1 &\mbox{ in }\O\label{CIdeltaborne}
\end{eqnarray}
and
\begin{eqnarray}
&0<u_m\leq c_\d \leq 1-\d<1 &\mbox{ in }Q_T.\label{Cdeltaborne}
\end{eqnarray}
Adding up (\ref{1d}) and (\ref{2d}) we deduce the equation
\begin{equation}\label{7d}
-div\bigg(\Mmu(u)\nabla p +\Frac{1}{\mu}
\ka(u)\nabla(\pc(u))\bigg)=\so_\d-\su_\d-\mint_\O(\so_\d-\su_\d
)dx.
\end{equation}
We formulate below an equivalent form of Problem $(S^{\mu}_\d)$.
To that purpose we define the global pressure, $\P$, by
$$\P := p + \calR(u) = p + \Int_0^u \frac{\ka(\tau)}{\ka(\tau)+\mu\kw(\tau)}\pc'(\tau)d\tau,$$
so that (\ref{7d}) gives
\begin{equation}\label{tn1d}
-div\bigg(\Mmu(u)\nabla \P
\bigg)=\so_\d-\su_\d-\mint_\O(\so_\d-\su_\d )dx.
\end{equation}
We rewrite the equation (\ref{1d}) of Problem $(S^{\mu}_\d)$ as
\begin{equation}\label{tn2d}
u_t =\Delta\psimu(u)+div\bigg(\fmu(u)\Mmu(u)\nabla \P\bigg)
+\fmu(c_\d)\so_\d - \fmu(u)\bigg(\su_\d+\mint_\O(\so_\d-\su_\d
)dx\bigg),
\end{equation}
where
\begin{equation}\label{phi}
\psimu(s)=-\Frac{1}{\mu}\Int_0^s\Frac{\ka(\tau)\kw(\tau)}{\Mmu(\tau)}\pc'(\tau)d\tau
\end{equation}
is continuous on $[0,1]$ and differentiable on $[0,1)$.
Multiplying (\ref{tn1d}) by $f^\mu(\umu_\d)$ and adding the result
to (\ref{tn2d}) we deduce that
$$
u_t =\Delta\psimu(u)+\Mmu(u)\nabla\fmu(u).\nabla \P+\bigg(\fmu(c_\d)-\fmu(u)\bigg)\so_\d.
$$
This yields a problem equivalent to $(S^{\mu}_\d)$, namely
$$
\initia \mbox{$(\tilde S^{\mu}_{\d})$} \left\{ \eqalijna{ & u_t
=\Delta\psimu(u)+\Mmu(u)\nabla\fmu(u).\nabla \P+[\fmu(c_\d)-\fmu(u)]\so_\d,
&\mbox{ in }\Omega \times (0,T),\cr &
 - div\bigg(\Mmu(u)\nabla\P\bigg)
 = \so_\d -\su_\d -\mint_\O(\so_\d-\su_\d )dx,&\mbox{ in
}\Omega \times (0,T), \cr
&\Int_\O\P(x,t)dx=\Int_\O\calR(u(x,t))dx,&\mbox{ for }t\in(0,T),\cr
 &\nabla \P.n=0, &\mbox{ on }\partial \Omega
\times (0,T),\cr &\nabla \psimu(u).n=0, &\mbox{ on }\partial \Omega
\times (0,T),\cr &u(x,0)=u_{0}^\d(x),& \mbox{ for }x \in \Omega.\cr
} \right. \eqno\eqalijnb{\latexeqno{t1d}\cr
              \latexeqno{t2d} \cr
               \latexeqno{t3d} \cr
                \latexeqno{t4d} \cr \latexeqno{t5d} \cr\latexeqno{tCId} \cr
                 }
$$
\noindent
In order to prove the existence of a smooth solution of
$(S^{\mu}_{\d})$, we introduce the set
$${\cal{K}}:=\{u\in C^{1+\a,\frac{1+\a}{2}}(\overline Q_T),~~u_m\leq u\leq 1-\d\},$$
where $\al\in(0,1)$, and we prove the following result
\begin{lemma}\label{leexistencetildeeps} Assume $(H_1)-(H_9)$ then there exists
$(\umu_{\d},\Pmu_{\d})$ solution of $(\tilde S^{\mu}_{\d})$ such
that
$$\umu_{\d}\in C^{2+\a,\frac{2+\a}{2}}(\overline Q_T),~~u_m\leq \umu_{\d}\leq 1-\d$$ and
$\Pmu_{\d}$, $\nabla \Pmu_{\d}$ $\in C^{1+\a,\frac{1+\a}{2}}(Q_T)$ and $\Delta\Pmu_{\d} \in
C^{\a,\frac{\a}{2}}(Q_T)$.
\end{lemma}
\proof Let $T^1$ be the map defined for all $V\in {\cal{K}}$ by $T^1(V)=W$, where $W$ is the unique solution of the
elliptic problem,
$$
(Q^{1}_V)\left\{
\begin{array}{lll}
& -div\bigg(\Mmu(V)\nabla W\bigg)=\so_\d-\su_\d-\mint_\O(\so_\d-\su_\d )dx,~~
&\mbox{ in }\Omega \times (0,T),\\
&\Int_\O W(x,t)dx=\Int_\O\calR(V(x,t))dx,&\mbox{ for }t\in (0,T),\\
&\nabla W.n=0, &\mbox{ on }\partial \Omega \times (0,T).
\end{array}\right.$$
By standard theory of elliptic system (see \cite{Lellip} Theorem 3.2 p 137), we have that
\begin{equation}\label{n8e}
|W|^{1+\a,\frac{1+\a}{2}}_{Q_T}+|\nabla
W|^{1+\a,\frac{1+\a}{2}}_{Q_T}\leq
D_1|V|^{1+\a,\frac{1+\a}{2}}_{Q_T}+D_2.
\end{equation}
For $W$ solution of $(Q^{1}_V)$ we consider $T^2$ defined by $T^2(W)=\hat V$, where $\hat V$ is the solution of the
parabolic problem,
$$
(Q^{2}_W)\left\{
\begin{array}{lll}
&\hat V_t=\Delta\psimu_\e(\hat V)+\Mmu(\hat
V)\nabla\fmu(\hat V).\nabla W +[\fmu(c_\d)-\fmu(\hat V)]\so_\d,
&\mbox{ in }\Omega \times (0,T),\\
&\nabla \psimu_\e(\hat V).n=0, &\mbox{ on }\partial \Omega \times (0,T),\\
&\hat V(x,0)=u_0^\d(x),&\mbox{ in }\Omega. \end{array}\right.$$
\noindent
>From the standard theory of parabolic equations, we have that
\begin{equation}\label{n8f}
|\hat V|^{2+\a,\frac{2+\a}{2}}_{Q_T}\leq
D_3\bigg(|W|^{1+\a,\frac{1+\a}{2}}_{Q_T}+|\nabla W|^{1+\a,\frac{1+\a}{2}}_{Q_T}\bigg) + D_4.\label{n8ebis}\\
\end{equation}
Moreover defining by $\cal L$ the parabolic operator arising in
$(Q^2_W)$, namely
$${\cal {L}}(\hat V)(x,t):=\hat V_t-\Delta\psimu_\e(\hat V)
-\Mmu(\hat V)\nabla\fmu(\hat V).\nabla W
-[\fmu(c_\d)-\fmu(\hat V)]\so_\d,$$ we remark that
(\ref{fmucroissante}), the property (\ref{Cdeltaborne}) of
$c_\d$ and the fact that, by (\ref{sdeltaborne}), $\so_\d$ is positive imply that
\begin{equation}\label{calO}
{\cal {L}}(u_m)\leq 0\mbox{ and }{\cal {L}}(1-\d)\geq 0.
\end{equation}
Setting $T:=T^2\circ T^1$, the inequalities (\ref{calO}) ensure that $T$ maps the convex set ${\cal{K}}$ into
itself. Moreover we deduce from (\ref{n8f}) that $T({\cal{K}})$ is relatively compact in ${\cal{K}}$.
\\Next, we check that $T$ is continuous. Suppose that a sequence $(V_m)_{m\in\N}$ converges to a limit $V\in
{\cal{K}}$ in $C^{1+\a,\frac{1+\a}{2}}(\overline Q_T)$, as $m\rightarrow\infty$. Since $(V_m)_{m\in\N}$ is bounded
in $C^{1+\a,\frac{1+\a}{2}}(\overline Q_T)$, it follows from (\ref{n8e}) that the sequence
$(W_m:=T^1(V_m))_{m\in\N}$, where $W_m$ is the solution of $(Q^1_{V_m})$, is bounded in
$C^{1+\a,\frac{1+\a}{2}}(\overline Q_T)$, so that as $m\rightarrow\infty$, $W_m$ converges to the unique solution
$W$ of Problem $(Q^1_V)$ in $C^{1+\beta,\frac{1+\beta}{2}}(\overline Q_T)$ for all $\beta\in(0,\a)$. Moreover $W\in
C^{1+\alpha,\frac{1+\alpha}{2}}(\overline Q_T)$. Further it also follows from (\ref{n8e}) that $(\nabla
W_m)_{m\in\N}$ is bounded in $C^{1+\a,\frac{1+\a}{2}}(\overline Q_T)$, so that the solution $\hat V_m=T^2(W_m)$ of
Problem $(Q^2_{W_m})$ is bounded in $C^{2+\a,\frac{2+\a}{2}}(\overline Q_T)$. Since $\hat V_m=T^2(W_m)=T(V_m)$,
$(T(V_m))_{m\in\N}$ converges to the unique solution $\hat V$ of Problem $(Q^2_{W})$ in
$C^{2+\beta,\frac{2+\beta}{2}}(\overline Q_T)$ for all $\beta\in(0,\a)$, as $m\rightarrow\infty$, so that $\hat
V=T^2(W)=T^2\circ T^1(V)$. Therefore we have just proved that $(T^2\circ T^1(V_m))_{m\in\N}$ converges to $T^2\circ
T^1(V)$ in $C^{2+\beta,\frac{2+\beta}{2}}(\overline Q_T)$ for all $\beta\in(0,\a)$, as $m\rightarrow\infty$, which
ensures the continuity of the map $T$. It follows from the Schauder fixed point theorem that there exists a solution
$(\umu_\d,\Pmu_\d)$ of $(\tilde{S}^{\mu}_{\d})$ such that
$$\umu_{\d}\in K\cap C^{2+\alpha,\frac{2+\alpha}{2}}(\overline Q_T) \mbox{ and $\Pmu_{\d}$, $\nabla \Pmu_{\d}$, $\in C^{1+\a,\frac{1+\a}{2}}(Q_T)$,
$\Delta \Pmu_{\d} \in C^{\a,\frac{\a}{2}}(Q_T)$.}$$ $\\
$ This concludes the proof of Lemma \ref{leexistencetildeeps}.
Moreover we deduce from Lemma \ref{leexistencetildeeps} the
existence of a solution of $(S^{\mu}_{\d})$, namely

\begin{corollary}\label{leexistenceeps}
Assume the hypotheses $(H_1)-(H_9)$ then there exists $(\umu_{\d},\pmu_{\d})$ solution of $(S^{\mu}_{\d})$ such that
$\umu_{\d}\in C^{2+\a,\frac{2+\a}{2}}(\overline Q_T)$,
\begin{equation}\label{borneu}
u_m\leq \umu_{\d}(x,t)\leq 1-\d
\end{equation}
and $\pmu_{\d}$, $\nabla \pmu_{\d}$ $\in C^{1+\a,\frac{1+\a}{2}}(Q_T)$, $\Delta \pmu_{\d} \in
C^{\a,\frac{\a}{2}}(Q_T)$.
\end{corollary}

\setcounter{equation}{0}
\section{A priori Estimates}
In view of (\ref{fmucroissante}) and (\ref{borneu}) we deduce the following bounds
\begin{eqnarray}
&&0=\fmu(0)\leq \fmu(\umu_{\d}(x,t))\leq 1=\fmu(1),\label{Eb1}\\
&&0=\fmu(0)\leq \fmu(c_{\d}(x,t))\leq 1=\fmu(1),\label{Eb1bis}\\
&&0<\kw(u_m)\leq \kw(\umu_{\d}(x,t))\leq \kw(1)=1,\label{Eb2}\\
&& 0=\ka(1)\leq \ka(\umu_{\d}(x,t))\leq \ka(0)=1,\label{Eb2bis}\\
&&0<\kw(u_m)\leq \Mmu(\umu_{\d}),\label{Eb3}\\
&&\pc(1)\leq \pc(\umu_{\d}(x,t))\leq \pc(0),\label{Eb4}\\
&& \pc(1)-\pc(0)\leq \calR(\umu_{\d}(x,t))\leq 0,\label{Eb5}\\
&& \pc(1)-\pc(0)\leq \calQ(\umu_{\d}(x,t))\leq 0,\label{Eb6}
\end{eqnarray}
for all $(x,t)\in\Omega\times(0,T)$. Next we state some essential
a priori estimates.
\begin{lemma}\label{leest1} Let $(\umu_{\d},\pmu_{\d})$ be a solution of Problem $(S^\mu_\d)$. There exists a positive constant $C$, which only depends on
$\Omega$, $\kw$, $\ka$ and $T$ such that
\begin{equation}\label{Est1}
\Int_0^T\IntO \ka(\umu_{\d})|\nabla\pmu_{\d}+\nabla
\pc(\umu_{\d})|^2dxdt\leq C\mu,
\end{equation}
\begin{equation}\label{Est1.1}
\Int_0^T\IntO |\nabla\pmu_{\d}|^2dxdt\leq C,
\end{equation}
and
\begin{equation}\label{Est1.2}
0\leq -\Int_0^T\IntO \nabla
g(\umu_{\d}).\nabla\pc(\umu_{\d})dxdt\leq C,
\end{equation}
\begin{equation}\label{Est1.2bis}
\Int_0^T\IntO |\nabla\zeta(\umu_{\d})|^2dxdt\leq C,
\end{equation}
\begin{equation}\label{Est1.2ter}
\Int_0^T\IntO |\nabla g(\umu_{\d})|^2dxdt\leq C.
\end{equation}
\end{lemma}
\proof We first prove (\ref{Est1}). Multiplying (\ref{tn1d}) by $\P=\pmu_{\d}+\calR(\umu_{\d})$ and integrating the
result on $Q_T=\Omega\times(0,T)$ we obtain
\begin{eqnarray}
&\Int_{Q_T}\Mmu(\umu_{\d})|\nabla(\pmu_\d+\calR(\umu_\d))|^2\leq
\Frac{1}{h}\Int_\QT(\so_\d-\su_\d)^2
+h\Int_\QT(\pmu_{\d}+\calR(\umu_{\d}))^2,\label{E1}
\end{eqnarray}
for all $h>0$. Moreover we have by Poincar\'e-Wirtinger inequality
that
$$\Int_\QT(\pmu_{\d}+\calR(\umu_{\d}))^2\leq C_1\bigg[\Int_\QT|\nabla(\pmu_{\d}+\calR(\umu_{\d}))|^2+
\bigg(\Int_\QT\pmu_{\d}+\calR(\umu_{\d})\bigg)^2 \bigg].$$ Using
(\ref{3d}) and (\ref{Eb5}), it follows that
$$\Int_\QT(\pmu_{\d}+\calR(\umu_{\d}))^2\leq C_1\Int_\QT|\nabla(\pmu_{\d}+\calR(\umu_{\d}))|^2+ C_2,$$
which we substitute into (\ref{E1}) with $h=\Frac{\kw(u_m)}{2C_1}$ to deduce, also in view of (\ref{sdeltaborne})
and (\ref{Eb3}), that
\begin{equation}\label{E2}
\Int_{Q_T}|\nabla(\pmu_{\d}+\calR(\umu_{\d}))|^2\leq C_3\mbox{ and
}\Int_{Q_T}|\pmu_{\d}+\calR(\umu_{\d})|^2\leq C_3.
\end{equation}
Furthermore multiplying (\ref{1d}) by $\pmu_{\d}$ and (\ref{2d}) by $\pmu_{\d}+\pc(\umu_{\d})$, adding up both
results and integrating on $Q_T$ we obtain
\begin{eqnarray}
&-\Int_{Q_T}(\umu_{\d})_t\pc(\umu_{\d})+\Int_{Q_T}
\kw(\umu_{\d})|\nabla\pmu_{\d}|^2+\Frac{1}{\mu}\ka(\umu_{\d})\big|\nabla\pmu_{\d}+\nabla
\pc(\umu_{\d})\big|^2=I,~~~~\label{E3}
\end{eqnarray}
where
$$\begin{array}{ll}
I:=&\Int_{Q_T}\bigg[\fmu(c_\d)\so_\d-\fmu(\umu_{\d})\bigg(\su_\d +\mint_\O(\so_\d-\su_\d )\bigg)\bigg]\pmu_{\d} dxdt\\
&+\Int_{Q_T}\bigg[(1-\fmu(c_\d))\so_\d-\bigg(1-\fmu(\umu_{\d})\bigg)\bigg(\su_\d+\mint_\O(\so_\d-\su_\d )
\bigg)\bigg](\pmu_{\d}+\pc(\umu_{\d}))dxdt.
\end{array}
$$
We check below that first term on the left-hand-side of (\ref{E3}) and $I$ are bounded. Denoting by ${\cal{P}}_c$ a
primitive of $\pc$ we have that
$$\Int_{Q_T}\pc(\umu_{\d})(\umu_{\d})_t=\Int_\O\Int_0^T\Frac{\partial}{\partial
t}\big[{\cal{P}}_c(\umu_{\d})\big].$$ Since ${\cal{P}}_c$ is continuous and $\umu_{\d}$ is bounded this gives
\begin{equation}\label{E4}
\bigg|\Int_{Q_T}\pc(\umu_{\d})(\umu_{\d})_tdxdt\bigg|\leq C_4.
\end{equation}
Moreover we have using (\ref{3d}) and (\ref{P+Q}) that
\begin{eqnarray}
&I=&\Int_{Q_T}\bigg(\pmu_{\d}+\calR(\umu_{\d})\bigg)(\so_\d-\su_\d)
dxdt\\
&&-\Int_{Q_T}\calR(\umu_{\d})\bigg[\fmu(c_\d)\so_\d-\fmu(\umu_{\d})\su_\d+\bigg(1-\fmu(\umu_\d)\bigg)\bigg( \mint_\O(\so_\d-\su_\d) \bigg)\bigg] dxdt\nonumber\\
&&+\Int_{Q_T}\bigg[(1-\fmu(c_\d))\so_\d-\bigg(1-\fmu(\umu_{\d})\bigg)\bigg(\su_\d
+\mint_\O(\so_\d-\su_\d) \bigg)
\bigg]\bigg[\calQ(\umu_{\d})+\pc(0)\bigg] dxdt.\nonumber
\end{eqnarray}
In view of $(H_5)$, (\ref{sdeltaborne}), (\ref{Eb1}),
(\ref{Eb1bis}), (\ref{Eb5}) and (\ref{Eb6}) we obtain
$$
I\leq C_5\Int_{Q_T}|\pmu_{\d}+\calR(\umu_{\d})|^2+ C_6.
$$
This together with (\ref{E2}) yields $I\leq C_5C_3+C_6$. Substituting this into (\ref{E3}) and also using (\ref{E4})
we obtain that
\begin{equation}\label{E5}
\Int_{Q_T} \kw(\umu_{\d})|\nabla\pmu_{\d}|^2+\Frac{1}{\mu}\ka(\umu_{\d})\big|\nabla\pmu_{\d}+\nabla
\pc(\umu_{\d})\big|^2dxdt\leq C_7,
\end{equation}
which implies (\ref{Est1}).
In view of (\ref{Eb2}), we also deduce from (\ref{E5}) the estimate (\ref{Est1.1}). \\
Next we prove (\ref{Est1.2}). By the definition (\ref{defg}) of $g$, we obtain from (\ref{7d}) that
\begin{equation}\label{E5b}
-div\bigg(\Mmu(\umu_{\d})\nabla\pmu_{\d}\bigg)
+\Frac{1}{\mu}\Delta g(\umu_{\d})=\so_\d-\su_\d-\mint_\O(\so_\d-\su_\d)dx.
\end{equation}
Multiplying (\ref{E5b}) by $\fmu(\umu_{\d})$ and subtracting the result from $(\ref{1d})$ we deduce that
\begin{equation}\label{E5t}
(\umu_{\d})_t=\Frac{1}{\mu}\fmu(\umu_{\d})\Delta
g(\umu_{\d})+div\bigg(k_w(\umu_{\d})\nabla \pmu_{\d}\bigg)
-\fmu(\umu_{\d})div\bigg(\Mmu(\umu_{\d})\nabla(\pmu)\bigg)+ \so_\d\big[\fmu(c_\d)-\fmu(\umu_{\d})\big].
\end{equation}
Moreover using the definition (\ref{deffmu}) of $\fmu$ and $\Mmu$ we note that
$$
\begin{array}{lll}
&div\bigg(\Mmu(\umu_{\d})\fmu(\umu_{\d})\nabla \pmu_{\d}\bigg)&=
div\bigg(k_w(\umu_{\d})\nabla \pmu_{\d}\bigg)\\
&&=\Mmu(\umu_{\d})\nabla(\fmu(\umu_{\d})).\nabla\pmu_{\d}
+\fmu(\umu_{\d})div\bigg(\Mmu(\umu_{\d})\nabla(\pmu)\bigg),
\end{array}
$$
which we substitute into (\ref{E5t}) to obtain
\begin{eqnarray}
(\umu_{\d})_t-\Frac{1}{\mu}\fmu(\umu_{\d})\Delta
g(\umu_{\d})- \Mmu(\umu_{\d})\nabla(\fmu(\umu_{\d})).\nabla\pmu_{\d}
=\so_\d\big[\fmu(c_\d)-\fmu(\umu_{\d})\big].\label{E6}
\end{eqnarray}
We set
\begin{equation}\label{defD}
D^\mu(a):=\pc(a)\fmu(a)-\calQ(a),
\end{equation}
for all $a\in[0,1]$, so that by the definition (\ref{Q}) of $\calQ$ we have $\nabla
D^\mu(\umu_{\d})=\pc(\umu_{\d})\nabla(\fmu(\umu_{\d}))$. Substituting this into (\ref{E6}), which we have multiplied
by $\pc(\umu_{\d})$, we deduce that
\begin{equation}\label{E7}
\pc(\umu_{\d})(\umu_{\d})_t -\Frac{1}{\mu}\fmu(\umu_{\d})\pc(\umu_{\d})\Delta
g(\umu_{\d})-\Mmu(\umu_{\d})\nabla
D^\mu(\umu_{\d}).\nabla\pmu_{\d}=\pc(\umu_{\d})\so_\d\big[\fmu(c_\d)-\fmu(\umu_{\d})\big].
\end{equation}
Multiplying (\ref{E5b}) by $D^\mu(\umu_{\d})$, adding the result
to (\ref{E7}) and also using the fact that
$$div\bigg(\Mmu(\umu_{\d})D^\mu(\umu_{\d})\nabla\pmu_{\d}\bigg)=
\Mmu(\umu_{\d})\nabla D^\mu(\umu_{\d}).\nabla\pmu_{\d}
+D^\mu(\umu_{\d})div\bigg(\Mmu(\umu_{\d})\nabla\pmu_{\d}\bigg),$$
we deduce that
\begin{eqnarray}
&\pc(\umu_{\d})(\umu_{\d})_t-\Frac{1}{\mu}\bigg(\fmu(\umu_{\d})\pc(\umu_{\d})-D^\mu(\umu_{\d})\bigg)\Delta
g(\umu_{\d})
-div\bigg(\Mmu(\umu_{\d})D^\mu(\umu_{\d})\nabla\pmu_{\d}\bigg)
\nonumber\\
&=\pc(\umu_{\d})\so_\d\bigg[\fmu(c_\d)-\fmu(\umu_{\d})\bigg]+
D^\mu(\umu_{\d})\bigg(\so_\d-\su_\d-\mint_\O(\so_\d-\su_\d)\bigg). \label{E8}\end{eqnarray} Integrating (\ref{E8})
on $Q_T$ and using the fact that the definition (\ref{defD}) of $D^\mu$ implies
$$\pc(\umu_{\d})\fmu(\umu_{\d})-D^\mu(\umu_{\d})=\calQ(\umu_{\d}),$$
we obtain
\begin{eqnarray}
\Int_{Q_T}\pc(\umu_{\d})(\umu_{\d})_tdxdt
-\Frac{1}{\mu}\Int_{Q_T}\calQ(\umu_{\d})\Delta
g(\umu_{\d})dxdt=J,\label{E9}
\end{eqnarray}
where
$$
\begin{array}{ll}
J:=&\Int_{Q_T}\pc(\umu_{\d})\so_\d\big[\fmu(c_\d)-\fmu(\umu_{\d})\big]dxdt\\
&+\Int_{Q_T}\bigg(\pc(\umu_{\d})\fmu(\umu_{\d})-\calQ(\umu_{\d})\bigg)
\bigg(\so_\d-\su_\d-\mint_\O(\so_\d-\su_\d)\bigg)dxdt.
\end{array}
$$
It follows from (\ref{Eb1}), (\ref{Eb1bis}), (\ref{Eb4}), (\ref{Eb6}) and
(\ref{sdeltaborne}) that $|J|\leq C_8$. Substituting this
into (\ref{E9}) and also using (\ref{E4}) we obtain that
\begin{equation}\label{E11}
0\leq -\Frac{1}{\mu}\Int_{Q_T}\nabla\calQ(\umu_{\d}).\nabla(
g(\umu_{\d}))dxdt\leq C_9.
\end{equation}
Furthermore we remark that
$$
\Frac{1}{\mu}\fmu(\umu_{\d})\geq \Frac{\kw(u_m)}{2},
$$
which together with (\ref{E11}) and the fact that $\nabla \calQ(\umu_{\d})=\fmu(\umu_{\d})\nabla\pc(\umu_{\d})$
yields
\begin{equation}\label{E12}
0\leq -\Int_{Q_T}\nabla\pc(\umu_{\d})\nabla(g(\umu_{\d}))dxdt\leq C_{10}.
\end{equation}
By the definition (\ref{defzeta}) of $\zeta$, we have $-\nabla\pc(\umu_{\d})\nabla
g(\umu_{\d})=|\nabla\zeta(\umu_{\d})|^2$. This together with (\ref{E12}) implies (\ref{Est1.2}) and
(\ref{Est1.2bis}), which in view of (\ref{Eb2bis}) gives (\ref{Est1.2ter}). This completes the proof of Lemma
\ref{leest1}. $\\ $$\\ $ In what follows we give estimates of differences of space translates of $\pmu_{\d}$ and
$g(\umu_{\d})$. We set for $r\in\R^+$ sufficiently small:
$$\Omega_r=\{x\in\Omega,~~B(x,2r)\subset\Omega\}.$$
\begin{lemma}\label{leest2} Let $(\umu_{\d},\pmu_{\d})$ be a solution of Problem $(S^\mu_\d)$;
there exists a positive constant $C$ such that
\begin{equation}\label{Est2}
\Int_0^T\Int_{\O_r}\bigg|\pmu_{\d}(x+\xi,t)-\pmu_{\d}(x,t)\bigg|^2(x,t)dxdt\leq
C\xi^2
\end{equation}
and
\begin{equation}\label{Est2.2}
\Int_0^T\Int_{\O_r}\bigg|g(\umu_{\d})(x+\xi,t)-g(\umu_{\d})(x,t)\bigg|^2dxdt\leq
C\xi^2,
\end{equation}
where $\xi\in\R^N$ and $|\xi|\leq 2r.$
\end{lemma}
\proof The inequalities (\ref{Est2}) and (\ref{Est2.2}) follow from (\ref{Est1.1}) and (\ref{Est1.2ter})
respectively. $\\ $$\\ $ Next we estimate differences of time translates of $g(\umu_{\d})$.
\begin{lemma}\label{leest3} Let $(\umu_{\d},\pmu_{\d})$ be a solution of Problem $(S^\mu_\d)$ then
there exists a positive constant $C$ such that
\begin{equation}\label{Est3}
\Int_0^{T-\tau}\Int_{\O}\big[g(\umu_{\d})(x,t+\tau)-g(\umu_{\d})(x,t)\big]^2dxdt\leq
C\tau,
\end{equation}
for all $\tau\in(0,T)$.
\end{lemma}
\proof We set
$$A(t):=\Int_\O[g(\umu_{\d})(x,t+\tau)-g(\umu_{\d})(x,t)]^2dx.$$
Since $g$ is a non decreasing Lipschitz continuous function with the Lipschitz constant $C_g$ we have that
\begin{eqnarray}
A(t)&\leq&  C_g\Int_\O[g(\umu_{\d}(x,t+\tau))-g(\umu_{\d}(x,t))][\umu_{\d}(x,t+\tau)-\umu_{\d}(x,t)]dx\nonumber\\
&\leq & C_g\Int_\O[g(\umu_{\d}(x,t+\tau))-g(\umu_{\d}(x,t))]\left[\Int_t^{t+\tau}(\umu_\d)_t(x,\theta)d\theta\right]dx\nonumber\\
&\leq&  C_g\Int_\O\Int_t^{t+\tau}\bigg[g(\umu_{\d}(x,t+\tau))-g(\umu_{\d}(x,t))\bigg]\nonumber\\
&~~~&~\bigg[div(\kw(\umu_{\d})\nabla\pmu_{\d})+\fmu(c_\d)\so_\d-\fmu(\umu_\d)
\bigg(\su_\d+\mint_\O(\so_\d(y)-\su_\d(y))dy\bigg)\bigg](x,\theta) d\theta dx,\nonumber
\end{eqnarray}
where we have used (\ref{1d}). Integrating by parts this gives
\begin{eqnarray}
&A(t)\leq
C_g\bigg\{\Int_t^{t+\tau}\Int_\O\bigg|\kw(\umu_{\d})(x,\theta)\nabla\pmu_{\d}(x,\theta)
\nabla g(\umu_{\d})(x,t+\tau)\bigg|dxd\theta \nonumber\\
&+\Int_t^{t+\tau}\Int_\O\bigg|\kw(\umu_{\d})(x,\theta)\nabla\pmu_{\d}(x,\theta)\nabla
g(\umu_{\d})(x,t)\bigg|
dxd\theta \nonumber\\
&+\bigg|\Int_\O\bigg[g(\umu_{\d})(x,t+\tau)-g(\umu_{\d})(x,t)\bigg]K(x,t,\tau)dx
\bigg|\bigg\},\label{Est20}
\end{eqnarray}
where
\begin{equation}\label{defK}
K(x,t,\tau):=\Int_t^{t+\tau}\bigg(\fmu(c_\d(x,\theta))\so_\d(x)-
\fmu(u_\d(x,\theta))\bigg[\su_\d(x)+\mint_\O(\so_\d(y)-\su_\d(y))dy\bigg]\bigg)d\theta.
\end{equation}
Next we estimate the right hand side of (\ref{Est20}). Using (\ref{Eb2}) we have that
\begin{eqnarray}
&&\Int_t^{t+\tau}\Int_\O\bigg|\kw(\umu_{\d})(x,\theta)\nabla\pmu_{\d}(x,\theta)
\nabla g(\umu_{\d})(x,t+\tau)\bigg|dxd\theta\nonumber\\
&&\leq
\Frac{1}{2}\bigg(\Int_t^{t+\tau}\Int_\O|\nabla\pmu_{\d}(x,\theta)|^2dxd\theta+
\Int_t^{t+\tau}\Int_\O|\nabla
g(\umu_{\d})(x,t+\tau)|^2dxd\theta\bigg)\nonumber\\
&&\leq
\Frac{1}{2}\bigg(\Int_t^{t+\tau}\Int_\O|\nabla\pmu_{\d}(x,\theta)|^2dxd\theta+
\tau\Int_\O|\nabla g(\umu_{\d})(x,t+\tau)|^2dx\bigg).
\label{Est21}
\end{eqnarray}
Similarly we have that
\begin{eqnarray}
&&\Int_t^{t+\tau}\Int_\O\bigg|\kw(\umu_{\d})(x,\theta)\nabla\pmu_{\d}(x,\theta)\nabla
g(\umu_{\d})(x,t)\bigg|dxd\theta\nonumber\\
&&\leq \Frac{1}{2}\bigg(\Int_t^{t+\tau}\Int_\O|\nabla\pmu_{\d}(x,\theta)|^2dxd\theta+ \tau\Int_\O|\nabla
g(\umu_{\d})(x,t)|^2dx\bigg). \label{Est22}
\end{eqnarray}
Moreover using (\ref{Eb1}) and (\ref{Eb1bis}) we obtain from the definition (\ref{defK}) of $K$ that
$$|K(x,t,\tau)|\leq
\Int_t^{t+\tau}\bigg[|\so_\d|+|\su_\d|+\mint_\O|\so_\d-\su_\d|dx\bigg]d\theta\leq
\bigg[|\so_\d|+|\su_\d|+\mint_\O|\so_\d-\su_\d|dx\bigg]\tau.$$ This together with $(\ref{sdeltaborne})$ and the fact
that the function $g(\umu_{\d})$ is bounded uniformly on $\mu$ and $\d$ yields
\begin{equation}\label{Est22.2}
\bigg|\Int_\O\bigg[g(\umu_{\d})(x,t+\tau)-g(\umu_{\d})(x,t)\bigg]K(x,t,\tau)\bigg|dx\leq \tilde C\tau.
\end{equation}
Substituting (\ref{Est21}), (\ref{Est22}) and (\ref{Est22.2}) into (\ref{Est20}) we deduce that
\begin{eqnarray}
A(t)&\leq
&C_g\bigg(\Int_t^{t+\tau}\Int_\O|\nabla\pmu_{\d}(x,\theta)|^2dxd\theta+
\Frac{\tau}{2}\Int_\O|\nabla
g(\umu_{\d})(x,t+\tau)|^2dx\nonumber\\
&&+ \Frac{\tau}{2}\Int_\O|\nabla g(\umu_{\d})(x,t)|^2dx+\tilde C\tau\bigg), \nonumber
\end{eqnarray}
which we integrate on $[0,T-\tau]$ to obtain
\begin{eqnarray}
\Int_0^{T-\tau}A(t)dt&\leq &C_g \bigg(
\Int_0^{T-\tau}\Int_t^{t+\tau}\Int_\O|\nabla\pmu_{\d}(x,\theta)|^2dxd\theta
dt+
\tau \Int_0^{T}\Int_\O|\nabla g(\umu_{\d})|^2dxdt+\tilde C\tau T\bigg)\nonumber\\
&\leq &C_g \bigg( \tau\Int_0^{T}\Int_\O|\nabla\pmu_{\d}(x,\theta)|^2dxd\theta + \tau \Int_0^{T}\Int_\O|\nabla
g(\umu_{\d})|^2dxdt+\tilde C\tau T\bigg).\nonumber
\end{eqnarray}
In view of (\ref{Est1.1}) and (\ref{Est1.2ter}) we deduce (\ref{Est3}), which completes the proof of Lemma
\ref{leest3}.
\setcounter{equation}{0}
\section{Convergence as $\d\downarrow 0$.}\label{deltavers0}
Letting $\d$ tend to 0, we deduce from the estimates given in Lemmas \ref{leest1} and \ref{leest2} the existence of
a weak solution of Problem $(S^\mu)$. More precisely, we have the following result,
\begin{lemma}\label{ledelta=0}
There exists a weak solution $(\umu,\pmu)$ of Problem $(S^\mu)$, which satisfies
\begin{equation}\label{Est1delta=0} \Int_0^T\IntO \bigg[\ka(\umu)\bigg]^2\bigg[\nabla\pmu+\nabla
\pc(\umu)\bigg]^2dxdt\leq C\mu,
\end{equation}
\begin{equation}\label{Est1.1delta=0}
\Int_0^T\IntO |\nabla\pmu|^2dxdt\leq C,
\end{equation}
\begin{equation}\label{Est1.2delta=0}
\Int_0^T\IntO |\nabla g(\umu)|^2dxdt\leq C,
\end{equation}
\begin{equation}\label{Est2.2delta=0}
\Int_0^T\Int_{\O_r}\big[g(\umu)(x+\xi,t)-g(\umu)(x,t)\big]^2dxdt\leq
C\xi^2,
\end{equation}
where $\xi\in\R^N$ and $|\xi|\leq 2r$. Moreover the following estimate of differences of time translates holds
\begin{equation}\label{Est3delta=0}
\Int_0^{T-\tau}\Int_{\O}\big[g(\umu)(x,t+\tau)-g(\umu)(x,t)\big]^2dxdt\leq
C\tau,
\end{equation}
for all $\tau\in(0,T)$.
\end{lemma}
\proof We deduce from (\ref{Est1.1}), (\ref{Est2.2}) and (\ref{Est3}) that there exist functions $\hat g^\mu$ and
$p^\mu$ and a subsequence $((u^{\mu}_\dn,p^{\mu}_\dn))_{n\in N}$ of weak solutions of Problem $(S^{\mu}_{\dn})$ such
that
\begin{eqnarray}
&&(g(u^{\mu}_\dn))_{n\in N}\mbox{ tends to }\hat g^\mu\mbox{ strongly in }L^2(Q_T), ~~~~~~\label{borneudelta=0.a}
\\
&&(p^{\mu}_\dn)_{n\in N}\mbox{ tends to }p^\mu\mbox{ weakly in }L^2(0,T;H^1(\O)),\nonumber
\end{eqnarray}
as $\d_n$ tends to zero. Thus for a subsequence, which we denote again by $\delta_n$, we have that
\begin{equation}\label{borneudelta=0.b}
(g(u^{\mu}_\dn))_{n\in N}\mbox{ tends to }\hat g^\mu\mbox{ for almost }(x,t)\in Q_T.
\end{equation}
Using the fact that g is bijective we deduce that
\begin{equation}\label{borneudelta=0.c} (u^{\mu}_\dn)_{n\in N} \mbox{ tends to $u^\mu:=g^{-1}(\hat g^\mu)$ strongly in $L^2(Q_T)$ and almost everywhere in }Q_T,
\end{equation}
as $\d_n$ tends to zero. Moreover we have in view of (\ref{Est1.2ter}) and (\ref{borneudelta=0.a}) that $\nabla
g(\umu_{\delta_n})$ tends to $\nabla g(\umu)$ weakly in $L^2(Q_T)$ as $\delta_n\downarrow 0$, so that by the
definition (\ref{defg}) of $g$
\begin{equation}\label{cvfaible1}
\ka(\umu_{\delta_n})\nabla p_c(\umu_{\delta_n}) \mbox{ tends to $\ka(\umu)\nabla p_c(\umu)$ weakly in $L^2(Q_T)$ as
$\delta_n\downarrow 0$.}
\end{equation}
Letting $\d_n$ tend to 0 in (\ref{borneu}) we deduce that
\begin{equation}\label{borneudelta=0.1}
u_m\leq \umu(x,t)\leq 1. \end{equation} Moreover we deduce from
(\ref{3d}) that
\begin{equation}\label{delta=0.2}
\Int_\O\pmu(x,t)dx=0, \mbox{ for almost every } t \in (0,T).
\end{equation} Multiplying (\ref{1d}) by $\varphi\in{\cal{C}}$,
integrating by parts and letting $\d_n$ tend to 0 we obtain
\begin{eqnarray}
&\Int_0^T\IntO u^\mu \varphi_tdxdt =&\Int_0^T\IntO\kw(u^\mu)\nabla p^\mu.\nabla \varphi dxdt-\Int_0^T\IntO\bigg( f^\mu(c)\so- f^\mu(u^\mu)\su \bigg) \varphi dxdt\nonumber\\
&&-\IntO u_0(x)\varphi(x,0)dx,\label{delta=0.3}
\end{eqnarray}
where we have used that $u_{0}^\d$ tends to $u_0$ in $L^2(\Omega)$ and that $c_\d$, $\so_\d$ and $\su_\d$ tend
respectively to $c$, $\so$ and $\su$ in $L^2(Q_T)$ as $\d\downarrow 0$. Similarly, multiplying (\ref{2d}) by
$\varphi\in{\cal{C}}$, integrating by parts and letting $\d_n$ tend to 0 we deduce that
\begin{eqnarray}
&&\Int_0^T\IntO\bigg(1-\umu\bigg)\varphi_tdxdt=\Frac{1}{\mu}\Int_0^T\IntO\bigg(\ka(\umu)\nabla\pmu+\nabla g(\umu)\bigg).\nabla \varphi dx dt\nonumber\\
&&-\Int_0^T\IntO\bigg((1-\fmu(c))\so-(1-\fmu(\umu))\su\bigg) \varphi dx dt-\IntO
\bigg(1-u_0(x)\bigg)\varphi(x,0)dx,~~\label{delta=0.4}
\end{eqnarray}
which since $\nabla g(\umu)=\ka(\umu)\nabla \pc(\umu)$ coincides with (\ref{defsol2}). Next we prove
(\ref{Est1delta=0}). We first check that
\begin{equation}\label{cvfaible2}
\ka(\umu_{\delta_n})\nabla \pmu_{\delta_n}\mbox{ tends to }\ka(\umu)\nabla \pmu \mbox{ weakly in }L^2(Q_T),
\end{equation}
as $\d_n$ tends to 0. Let $\varphi\in L^2(Q_T)$, we have that
\begin{eqnarray}
\bigg|\Int_{Q_T} \bigg(\ka(\umu_{\delta_n})\nabla \pmu_{\delta_n}-\ka(\umu)\nabla \pmu\bigg)\varphi ~dxdt\bigg| \leq
|I^1_{\delta_n}|+|I^2_{\delta_n}|,\label{delta=0.4bis}
\end{eqnarray}
where
$$
I^1_{\delta_n}:= \Int_{Q_T} \bigg(\ka(\umu_{\delta_n})-\ka(\umu)\bigg) \nabla \pmu_{\delta_n}\varphi~dxdt
$$and
$$I^2_{\delta_n}=\Int_{Q_T}
\ka(\umu)\varphi\bigg(\nabla \pmu_{\delta_n} -\nabla \pmu\bigg)~dxdt.
$$
Using the fact that $\nabla \pmu_{\delta_n}$ converges to $\nabla \pmu$ weakly in $L^2(Q_T)$ as $\delta_n\downarrow
0$, we deduce, since $\ka(\umu)\varphi\in L^2(Q_T)$, that
\begin{equation}\label{delta=0.4ter}
\mbox{ $|I^2_{\delta_n}|$ tends to $0$ as $\delta_n\downarrow 0$.}
\end{equation}
Moreover we have by (\ref{Est1.1}) that
\begin{eqnarray}
|I^1_{\delta_n}|&&\leq \bigg(\Int_{Q_T} \bigg|\ka(\umu_{\delta_n})-\ka(\umu)\bigg|^2\varphi^2dxdt\bigg)^{1/2}
\bigg(\Int_{Q_T}|\nabla \pmu_{\delta_n}|^2dxdt\bigg)^{1/2} \nonumber\\&&\leq C\bigg(\Int_{Q_T}
\bigg|\ka(\umu_{\delta_n})-\ka(\umu)\bigg|^2\varphi^2dxdt\bigg)^{1/2}.\nonumber
\end{eqnarray}
Since $\bigg|\ka(\umu_{\delta_n})-\ka(\umu)\bigg|^2\varphi^2\leq 4\varphi^2$ and since $\ka(\umu_{\delta_n})$ tends
to $\ka(\umu)$ almost everywhere, we deduce from the Dominated Convergence Theorem that $I^1_{\delta_n}$ tends to
$0$ as $\delta_n\downarrow 0$. This with (\ref{delta=0.4ter}) implies (\ref{cvfaible2}), which with
(\ref{cvfaible1}) gives that
\begin{equation}\label{cvfaible3}
\ka(\umu_{\delta_n})\bigg[\nabla \pmu_{\delta_n}+\nabla p_c(\umu_{\delta_n})\bigg]\mbox{ tends to
}\ka(\umu)\bigg[\nabla \pmu+\nabla p_c(\umu)\bigg] \mbox{ weakly in }L^2(Q_T).
\end{equation}
The functional $v\mapsto \Int_{Q_T}v^2dxdt$ is convex and lower semi continuous from $L^2(Q_T)$ to $\overline R$
therefore it is also weakly l.s.c. (see \cite{B} Corollary III.8) and thus we deduce from (\ref{Eb2bis}),
(\ref{Est1}) and (\ref{cvfaible3}) that
$$
\begin{array}{lll}
\Int_{Q_T}\bigg[\ka(\umu)\bigg]^2\bigg[\nabla \pmu +\nabla p_c(\umu)\bigg]^2dxdt
&\leq& \liminf_{\delta_n\downarrow 0
}\Int_{Q_T}\bigg[\ka(\umu_{\delta_n})\bigg]^2\bigg[\nabla \pmu_{\delta_n} +\nabla
p_c(\umu_{\delta_n})\bigg]^2dxdt\\
&\leq& \liminf_{\delta_n\downarrow 0}\Int_{Q_T}\ka(\umu_{\delta_n})\bigg[\nabla \pmu_{\delta_n} +\nabla
p_c(\umu_{\delta_n})\bigg]^2dxdt \\
&\leq &C\mu,
\end{array}
$$ which coincides with (\ref{Est1delta=0}). Finally, we deduce
respectively from (\ref{Est1.1}), (\ref{Est1.2ter}), (\ref{Est2.2}) and (\ref{Est3}) the estimates
(\ref{Est1.1delta=0}), (\ref{Est1.2delta=0}), (\ref{Est2.2delta=0}) and (\ref{Est3delta=0}). This concludes the
proof of Lemma \ref{ledelta=0}.

\setcounter{equation}{0}
\section{Convergence as $\mu\downarrow 0$.}\label{muvers0}
The goal of this section is to prove Theorem \ref{thlim}. We first deduce from the estimates (\ref{Est1.1delta=0}),
(\ref{Est2.2delta=0}) and (\ref{Est3delta=0}) that there exists a couple of functions $(u,p)$ and a subsequence
$((u^{\mu_n},p^{\mu_n}))_{n\in N}$ such that $$\begin{array}{ll}
&(u^{\mu_n})_{n\in N}\mbox{ tends to }u\mbox{ strongly in }L^2(Q_T)\mbox{ and almost everywhere in }Q_T ,\\
&(p^{\mu_n})_{n\in N}\mbox{ tends to }p\mbox{ weakly in }L^2(0,T;H^1(\O)),\\
\end{array}$$
as $\mu_n$ tends to zero. Moreover since
$$0\leq f^{\mu_n}(u^{\mu_n})\leq 1,$$
there exists a function $\hat f\in L^2(Q_T)$ with $0\leq \hat f\leq 1$ and a subsequence
$(f^{\mu_{n_m}}(u^{\mu_{n_m}}))_{n_m\in N}$ of $(f^{\mu_{n}}(u^{\mu_{n}}))_{n\in N}$ such that
$(f^{\mu_{n_m}}(u^{\mu_{n_m}}))_{n_m\in N}$ tends to $\hat f$ weakly in $L^2(Q_T)$ as $\mu_{n_m}$ tends to zero.
Moreover we deduce respectively from (\ref{borneudelta=0.1}) and (\ref{delta=0.2}) that $0\leq u\leq 1$ and that
$$\Int_\Omega p(x,t)dx=0,\mbox{ for almost every }t\in(0,T),$$
which gives (\ref{lim2}). As it is done in Section \ref{muvers0} in the proof of (\ref{Est1delta=0}), one can first
check that
$$
\ka(u^{\mu_{n_m}})(\nabla p^{\mu_{n_m}} +\nabla p_c(u^{\mu_{n_m}}))\mbox{ tends to }\ka(u)(\nabla p +\nabla
p_c(u))\mbox{ weakly in }L^2(Q_T),
$$
as $\mu_{n_m}\downarrow 0$ and then deduce from (\ref{Est1delta=0}) the estimate (\ref{lim1}). Furthermore letting $\mu_{n_m}$ tends to zero into
(\ref{delta=0.3}) we obtain, since $\lim_{\mu_{n_m}\downarrow 0}f^{\mu_{n_m}}(s)=\chi(s)$ for all $s\in[0,1]$, that
$$
\begin{array}{ll}
\Int_0^T\IntO u \varphi_tdxdt =&\Int_0^T\IntO\kw(u)\nabla p.\nabla \varphi dxdt-\Int_0^T\IntO\bigg(\chi(c)\so-\hat f\su \bigg) \varphi dxdt\nonumber\\
&-\IntO u_0(x)\varphi(x,0)dx,
\end{array}
$$
which coincides with (\ref{sollim1}) and concludes the proof of Theorem \ref{thlim}.

\setcounter{equation}{0}
\section{Numerical simulations}
\subsection{The saturation equation and the numerical algorithm}
In this section we present numerical simulations in one space dimension. To that purpose we apply the finite volume
method, which we present below. To begin with, we rewrite the equations (\ref{n1}) and (\ref{n2}) in the case that
$\Omega=(0,1)$; this gives for $(x,t)\in(0,1)\times(0,T)$
\begin{eqnarray}
&&u^\mu_t =\partial_x\bigg(\kw(u^\mu)\partial_x
p^\mu\bigg)+\fmu(c)\so-\fmu(u^\mu)\su, ~~\label{D1} \\
&&(1-u^\mu)_t =
\partial_x\bigg(\Frac{1}{\mu}\ka(u^\mu)\partial_x(p^\mu+\pc(u^\mu))\bigg)+ (1-\fmu(c))\so-(1-\fmu(u^\mu))\su.~~\label{D2}
\end{eqnarray}
Adding up both equations and using the boundary conditions (\ref{n3}) and (\ref{n5}) we obtain
\begin{equation}\label{D3}
\partial_x p^\mu=-\Frac{\ka(u^\mu)}{\ka(u^\mu)+\mu\kw(u^\mu)}\partial_x(\pc(u^\mu)).
\end{equation}
Substituting (\ref{D3}) into (\ref{D1}) yields
\begin{equation}\label{D4}
u^\mu_t
=-\partial_x\bigg[\fmu(u^\mu)\Frac{\ka(u^\mu)}{\mu}\partial_x(\pc(u^\mu))\bigg]+\fmu(c)\so-\fmu(u^\mu)\su.
\end{equation}
Moreover we deduce from (\ref{D3}) and the definition (\ref{R}) of $\calR$ that $\partial_x
p^\mu=-\partial_x(\calR(u^\mu))$, so that in view of (\ref{n4}) we have
\begin{equation}\label{D5}
p^\mu(x,t)=-\calR(u^\mu)(x,t)+\Int_0^1\calR(u^\mu)(y,t)dy.
\end{equation}
In the sequel, we compare numerically the solution $u^\mu$ of (\ref{D4}) with the solution $u$ of the limit equation
in the case that $u<1$, namely
\begin{equation}\label{D5.b}
u_t =-\partial_x\bigg(k_w(u)\partial_x\pc(u)\bigg)+ \chi(c)\bar s.
\end{equation}
We discretize the time evolution equation (\ref{D4}) together with the initial condition and the homogeneous Neumann
boundary condition. The time explicit finite volume scheme is defined by the following equations in which $\k>0$ and
$\i>0$ denote respectively  the time and the space step. \\
(i) The discrete initial condition is given for $i\in\{0,...,[1/\i]\}$ by
\begin{equation}\label{Condinitialdiscre}
[U^\mu]_i^{0}=u^\mu(i\i,0).
\end{equation}
(ii) For $i\in\{0,...,[1/\i]\}$ and for $n\in\{0,...,[T/\k]\}$ the discrete equation is given by
\begin{eqnarray}
\Frac{1}{\k}\bigg([U^\mu]_i^{n+1}-[U^\mu]_i^{n}\bigg)&=&
[F^\mu]_{i+1}^n-[F^\mu]_i^n
+\fmu(C_i^n)\overline{S}_i^n-\fmu([U^\mu]_i^n)\underline{S}_i^n,\label{D6}
\end{eqnarray}
where
$$[F^\mu]_i^n=-\Frac{1}{\i}\bigg(p_c([U^\mu]_{i+1}^n)-p_c([U^\mu]_{i}^n\bigg)
\Frac{\kw([U^\mu]_{i+1}^n)\ka([U^\mu]_{i}^n)}{\mu\kw([U^\mu]_{i+1}^n)+\ka([U^\mu]_{i}^n)}.$$

\noindent(iii) For $n\in\{0,...,[T/\k]\}$ the discrete Neumann condition is defined  by
\begin{equation}\label{Condborddiscre}
[F^\mu]_0^{n}=0 \mbox{ and } [F^\mu]_{[1/\i]}^{n}=0.
\end{equation}
The numerical scheme (\ref{Condinitialdiscre})-(\ref{Condborddiscre}) allows to build an approximate solution,
$u_{\i,\k}:[0,1]\times[0,T]\rightarrow \R$ for all $i\in\{0,...,[1/\i]\}$ and all $n\in\{0,...,[T/\k]\}$, which is
given by
\begin{equation}\label{schema}
u_{\i,\k}(x,t)=u_i^n,\mbox{ for all }x\in(i\i,(i+1)\i]\mbox{ and
for all }t\in (n\k,(n+1)\k].
\end{equation}
In order to also compute the pressures, we propose the following
discrete equation corresponding to (\ref{D5})
\begin{equation}\label{D7}
[P^\mu]_i^{n}=-{\cal{R}}([U^\mu]_i^{n})+\i\Sigma_{j=1}^{[1/\i]}{\cal{R}}([U^\mu]_j^{n}).
\end{equation}
Finally, setting $p^\mu_g(x,t)=p^\mu(x,t)+p_c(u^\mu)(x,t)$ we
deduce that
\begin{equation}\label{D8}
([P_g^\mu])_i^{n}=-{\cal{R}}([U^\mu]_i^{n}+p_c([U^\mu]_i^{n})
+\i\Sigma_{j=1}^{[1/\i]}{\cal{R}}([U^\mu]_j^{n}),
\end{equation}
for all $i\in\{0,...,[1/\i]\}$ and all $n\in\{0,...,[T/\k]\}$.
Similarly we propose a finite volume scheme corresponding to the
equation (\ref{D5.b}), namely
\begin{eqnarray}
\Frac{1}{\k}\bigg(U_i^{n+1}-U_i^{n}\bigg)= F_{i+1}^n-F_i^n
+\chi(C_i^n)\overline{S}_i^n\label{D6},
\end{eqnarray}
where
$$F_i^n=-\Frac{1}{\i}\bigg(p_c(U_{i+1}^n)-p_c(U_{i}^n\bigg)
\kw(U_{i+1}^n),$$ for all
$(i,n)\in\{0,...,[1/\i]\}\times\{0,...,[T/\k]\}$.

\subsection{Numerical tests}
For the numerical computation we take $\mu=10^{-8}$, $p_c(z)=0,1\sqrt{1-z}$, $k_a(z)=(1-z)^2$, $k_w(z)=\sqrt z$ and
$\so(z)=\delta_0(z)$, $\su(z)=\delta_1(z)$, where $\delta_a$ is the Dirac function at the point $a$. Furthermore
$u^\mu$ is
given by the line with crosses, $p^\mu_g$ is given by the lines with diam and the limit function $u$ corresponds to the continuous line.\\
$\\ $ \underline{\bf{First test case}:} The case that $c=0,7$ and
$u_0=1$ on $[0,1]$. We obtain at $t=0,01$ the following pictures
\begin{center}
\includegraphics[width=0.5\hsize]{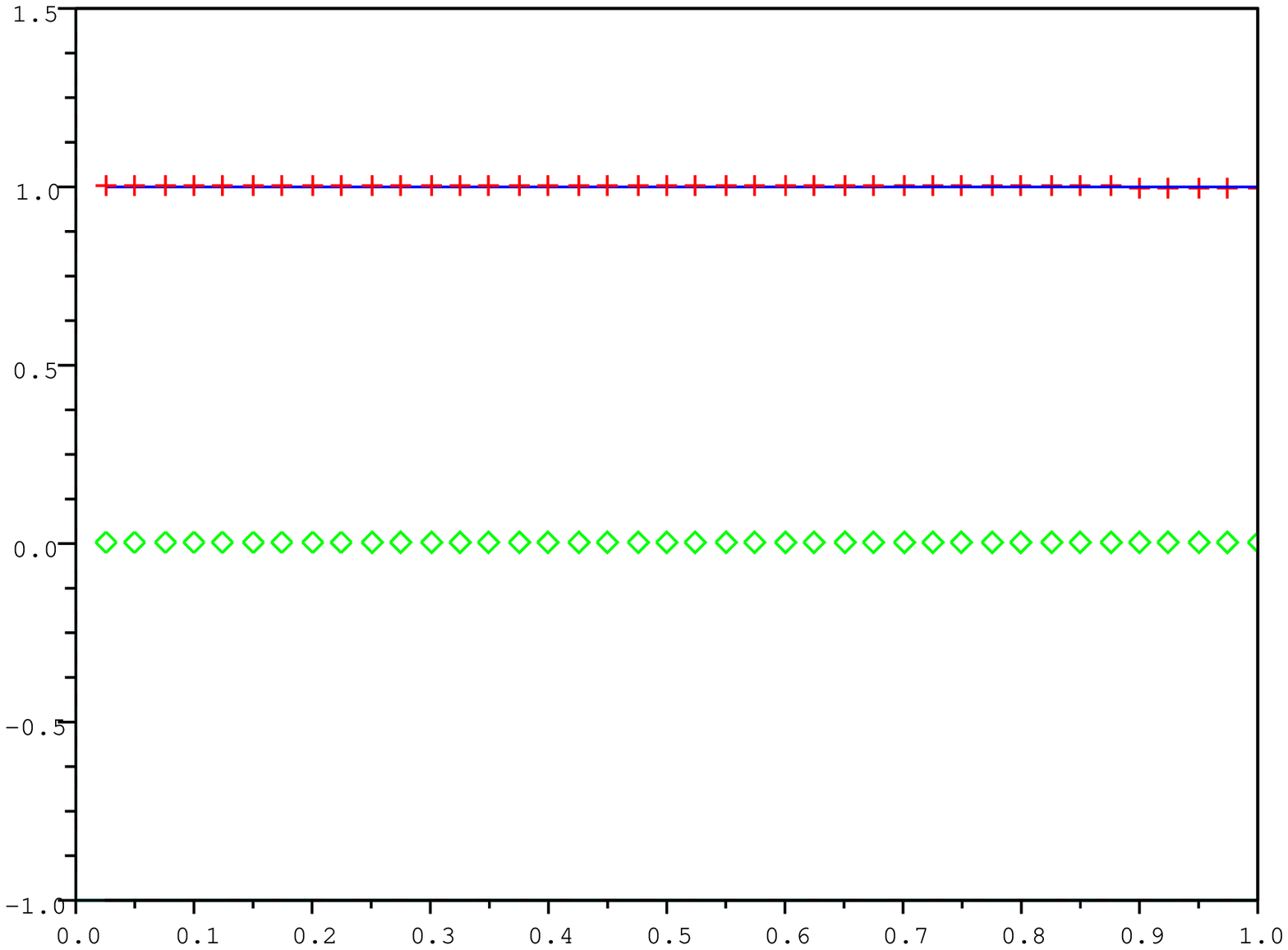}\\
\it{Figure 1 : t=0,01}
\end{center}
We note that, for $\mu$ small, the functions $u$ and $u^\mu$ are very close. Here we only start with water and
inject a mixture of water and air. The air immediately invades the whole domain. Figure 1 illustrates the result
which we proved in this paper, namely that $u^\mu$ tends to the solution $u$ of the limit equation (\ref{D5.b}) as
$\mu$ tends to 0 and moreover that the pressure
$p^\mu_a$ is constant. This is indeed the case since $u<1$.\\
$\\ $
\underline{\bf{Second test case}:} The case that $c=0,7$ and
$u_0(x)=\left\{\begin{array}{ll} 0,1\mbox{ on }[0,1/3]\\0,7\mbox{
on }(1/3,1]\end{array}\right.$. We obtain the following pictures
for $t=0,01$ and for $t=0,1$ respectively
\begin{center}
\includegraphics[scale=0.4]{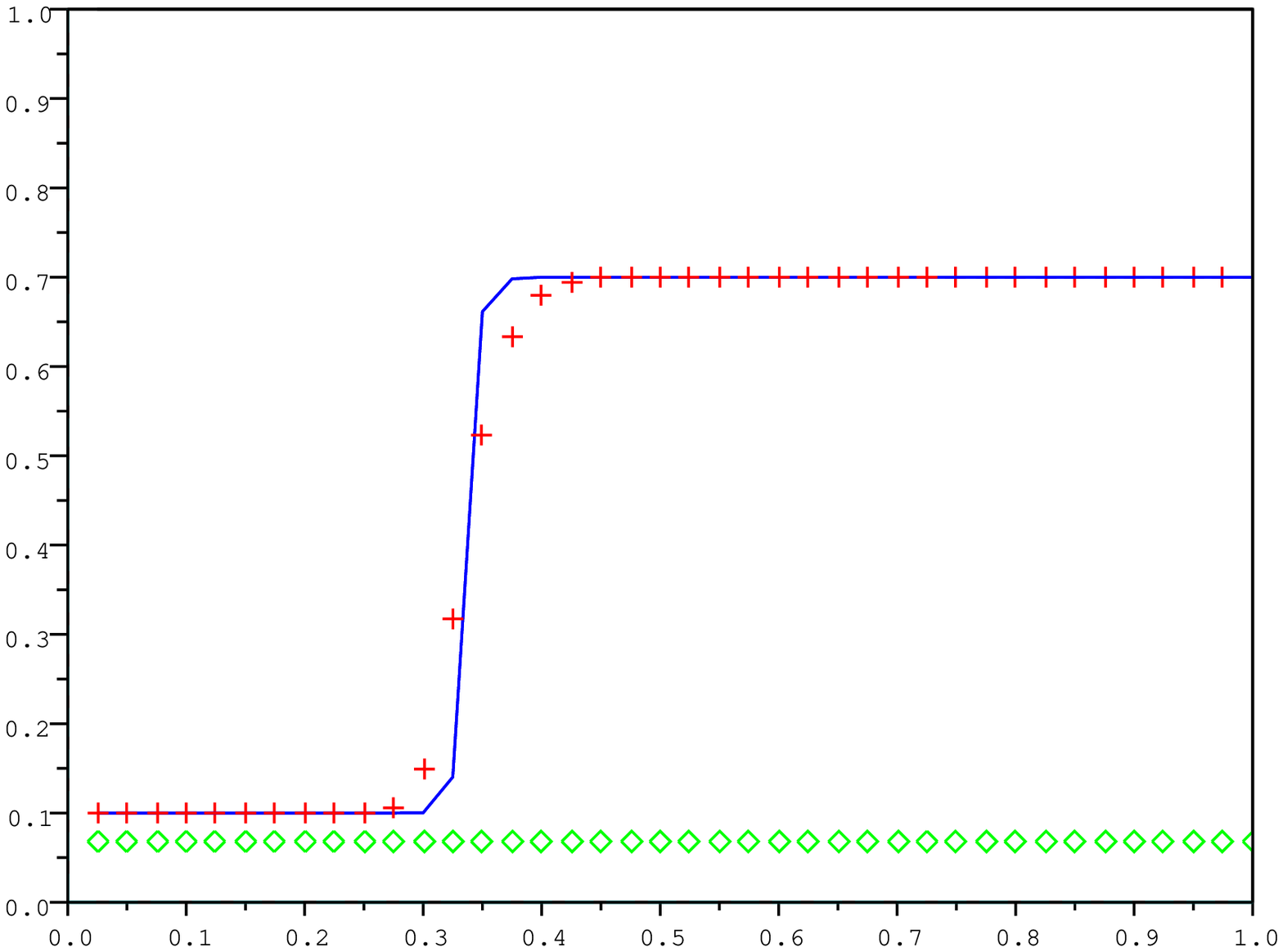}\\
\it{Figure 2 : t=0,01}
\end{center}
$\\ $
\begin{center}
\includegraphics[scale=0.4]{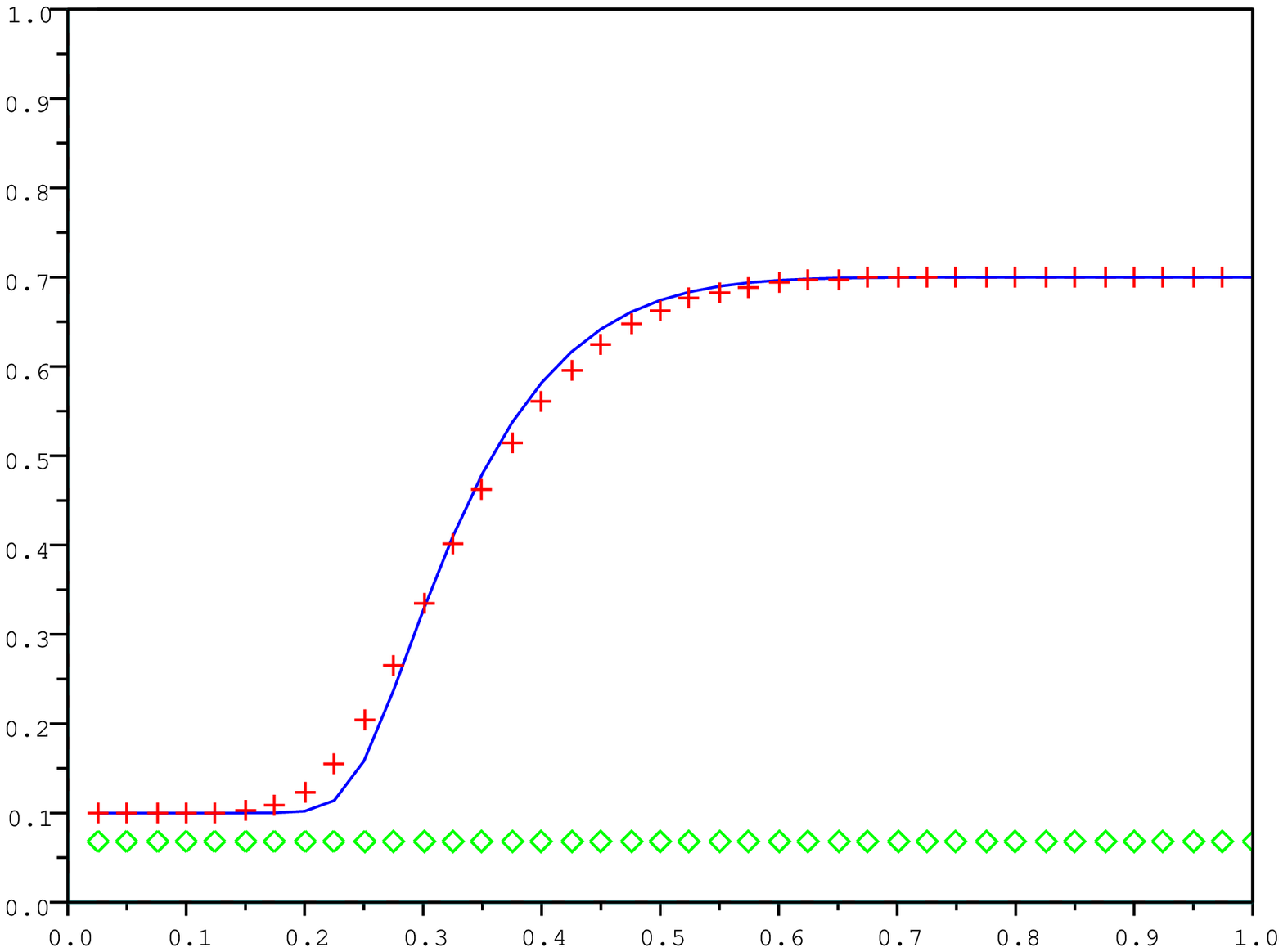}\\
\it{Figure 3 : t=0,1}
\end{center}
The injection of a mixture of water and air $(c=0.7)$ takes place in a region of low water saturation. We first
remark that both functions $u^\mu$ and $u$ evolve very slowly. Here again we have that $u(x,t)<1$ for all
$(x,t)\in(0,1)\times(0,T)$ and we remark that the graphs of the two functions $u^\mu$ and $u$ nearly
coincide.\\
$\\ $ \underline{\bf{Third test case}:} The case that $c=1$ and
$u_0(x)=\left\{\begin{array}{ll} 0,1\mbox{ on }[0,1/3]\\0,7\mbox{
on }(1/3,1]\end{array}\right.$. We obtain the following pictures
for $t=0,01$ and for $t=0,1$ respectively
\begin{center}
\includegraphics[scale=0.4]{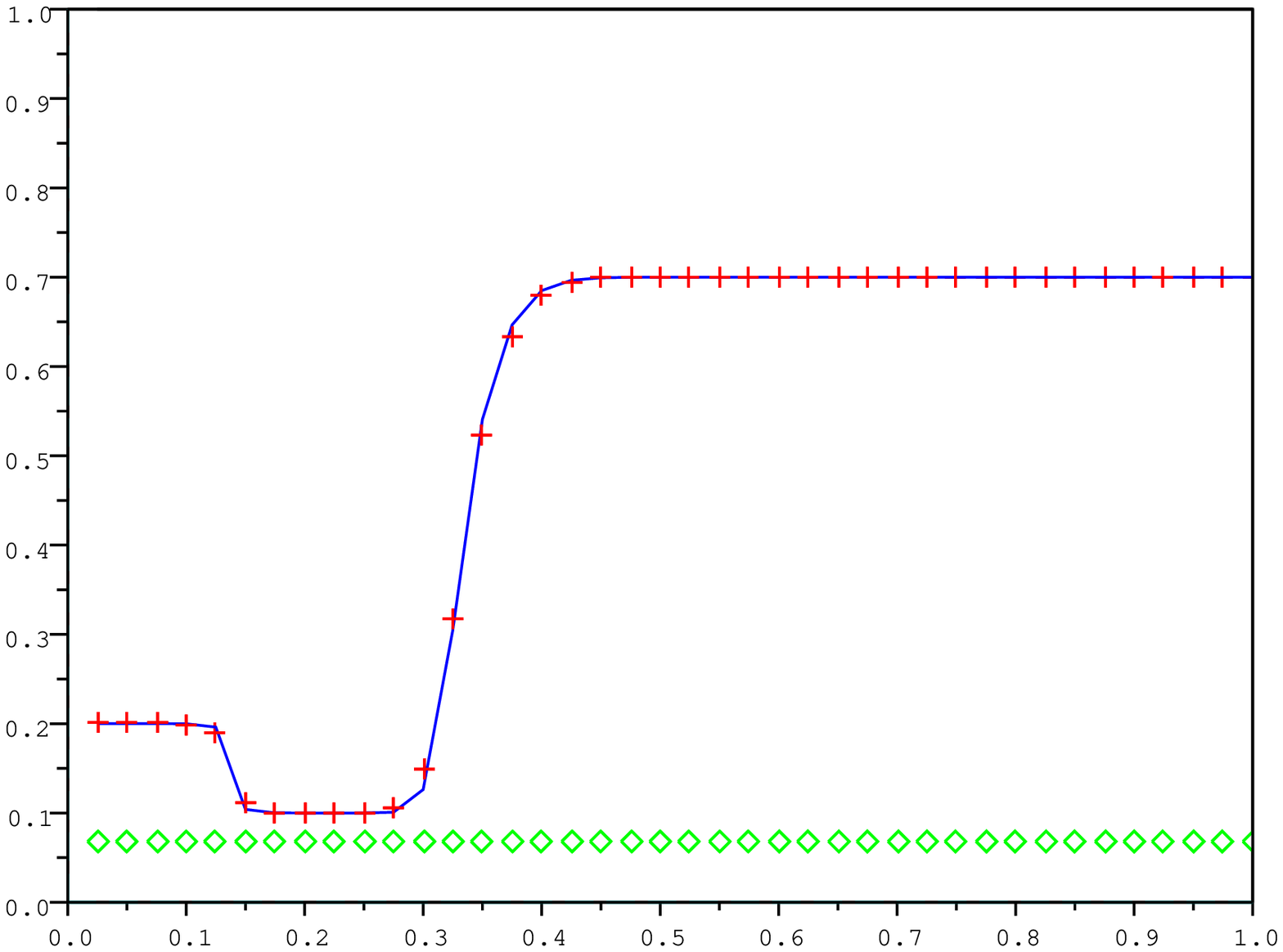}\\
\it{Figure 4 : t=0,01}
\end{center}
$\\ $
\begin{center}
\includegraphics[scale=0.4]{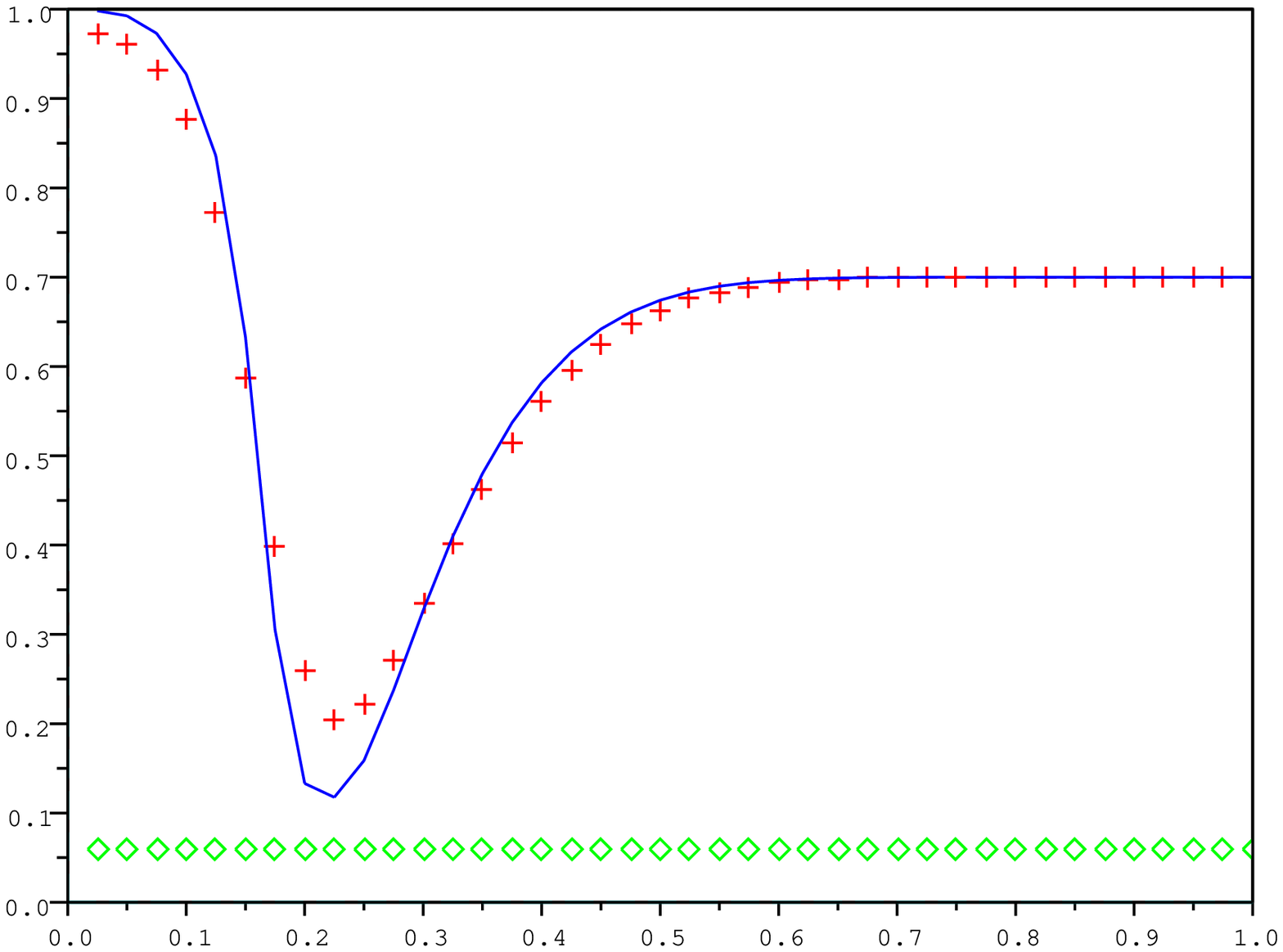}\\
\it{Figure 5 : t=0,1}
\end{center}
Here only water is injected; note that the saturation $\umu$ evolves rather fast.

\newpage \noindent
{\Large\bf Acknowledgement} \\
\vspace{0.05 cm} \\
We would like to thank Professor Rapha\`ele Herbin from the University of Provence for her interest in this work.

\end{document}